\documentstyle{amsppt}
\magnification=\magstep1
\hsize=6.5 true in
\vsize=9 true in
\topmatter
\title
Factoring Hecke polynomials modulo a prime\\
\endtitle
\author
J.B.~Conrey \\
D.W.~Farmer \\
P.J.~Wallace
\endauthor
\thanks Research of the first two authors supported 
by the American Institute of Mathematics.
Research of the first author supported in
part by a grant from the NSF.
A portion of this work resulted from the third author's 
Honors Thesis at Bucknell University.
\endthanks
\address
American Institute of Mathematics,
Palo Alto, CA 94306
\endaddress
\address
 Department of Mathematics, Oklahoma State University,
Stillwater, OK 74078
\endaddress
\address
 Department of Mathematics, Bucknell University,
Lewisburg, PA 17837
\endaddress
\abstract
Let $T_{p,k}^{N,\chi}(x)$ be 
the characteristic polynomial of the 
Hecke operator $T_p$ acting on the space of cusp forms $S_k(N,\chi)$.
We describe the factorization of $T_{p,k}^{N,\chi}(x)\bmod \ell$
as $k$ varies, and we explicitly calculate those factorizations
for $N=1$ and small $\ell$.
These factorizations are  used to deduce the
irreducibility of certain $T_{q,k}^{1,1}(x)$ from the
irreducibility of $T_{2,k}^{1,1}(x)$.
\endabstract
\endtopmatter
\document

\NoBlackBoxes

\def\mod{\bmod}

\def\Z {{\Bbb Z}}
\def\N {{\Bbb N}}

\def\Q{{\Bbb Q}}

\def\K{{\Bbb K}_k}

\def\({\left(}
\def\){\right)}

\def\newstrut{\hbox{\vrule height6.0pt depth 2.5pt width0pt}}
\def\newstrutA{\hbox{\vrule height10.0pt depth 3.5pt width0.4pt}}
\def\newstrutB{\hbox{\vrule height10.0pt depth 3.5pt width0.0pt}}
\def\tablewidthA{110pt}
\def\tablewidthB{230pt}
\def\tablewidthC{275pt}
\def\tablerule{\noalign{\hrule}}
\def\makeheadingfortableA{ &&\omit \hidewidth $\ell=5\ $\hidewidth &&
  \omit \hidewidth 0\hidewidth &&
  \omit \hidewidth 2\hidewidth&\cr\tablerule}
\def\makeheadingfortableB{ &&\omit \hidewidth $\ell=7\ $\hidewidth &&
  \omit \hidewidth 0\hidewidth &&
  \omit \hidewidth 2\hidewidth &&
  \omit \hidewidth 4\hidewidth\cr\tablerule}
\def\makeheadingfortableC{&&\omit \hidewidth $\ell=13\ $\hidewidth &&
  \omit \hidewidth ($a_1 \ldots a_{14}$)\hidewidth 
 %&& \omit \hidewidth \hidewidth
\cr\tablerule}

\def\TpkNchi{T_{p,k}^{N,\chi}(x)}
\head
1.  Introduction and statement of results
\endhead

Let
$S_k(N,\chi) $
be the space of holomorphic cusp forms of
integral weight $k$ for the 
Hecke congruence subgroup $\Gamma_0(N)$,
and denote by
$T_{p,k}^{N,\chi}(x)$ 
the characteristic polynomial of the 
Hecke operator $T_p$ acting on $S_k(N,\chi)$.
Let $S_k(N)=S_k(N,1)$ and set $T_{p,k}=T_{p,k}^{1,1}$.
A conjecture of Maeda 
asserts that $T_{p,k}(x)$ is irreducible and has
full Galois group over~$\Q$.
This conjecture is related to the nonvanishing of
$L$--functions \cite{KZ}\cite{CF}, and to
constructing base changes to totally
real number fields for level~1 eigenforms \cite{HM}.
Maeda's conjecture has been
checked for all primes $p<2000$ and weights $k\le 2000$
\cite{B}\cite{CF}\cite{FJ}.  

The methods which have been used to check Maeda's conjecture 
involve computing the factorization of $T_{2,k}(x) \bmod \ell$
for various $\ell$.  
One searches for enough %different 
factorizations to deduce
that $T_{2,k}(x)$ is irreducible and has full Galois
group over~$\Q$.  With a small amount of additional
calculation it is possible to deduce the same conclusion for other
$T_{p,k}(x)$.  The method of translating information about
$T_{2,k}(x)$ to information about $T_{p,k}(x)$ is related
to the following result of James and Ono.

\proclaim{Theorem\cite{JO}}  Suppose $q$ and $\ell$
are distinct primes.  Then 
% and $T_{p,k}(x)$ is irreducible mod~$\ell$.
$$
\#\left\{\hbox{prime } p<X\ |\ T_{p,k}^{N,\chi}(x)
\equiv T_{q,k}^{N,\chi}(x) \bmod \ell \right\}
\gg X/\log X.
$$
In particular, if  $\,T_{p,k}^{N,\chi}(x)$ 
is irreducible mod~$\ell$ for some~$p$, then 
the same holds for a positive proportion of primes~$p$.
\endproclaim

\break

The above result makes use of the Chebotarev density
theorem and deep results of Deligne,
Serre, and Shimura on the Galois representations
associated to modular forms.  
In this paper we describe how to determine, with a finite calculation, 
the factorization
of $T_{p,k}^{N,\chi}(x) \bmod \ell$ for any $k$.  We use
those factorizations to give an easy proof of the following
result, which is similar to the second statement in the theorem given above.

\proclaim{Theorem 1} If $\,T_{n,k}(x)$ is
irreducible and has full Galois group for some $n$,  
then $T_{p,k}(x)$ is
irreducible and has full Galois group for $p$
prime and
 $p\not \equiv \pm 1\bmod 5$, or
 $p\not \equiv \pm 1\bmod 7$.
\endproclaim

Note that the conclusion of Theorem~1 holds for $5/6$ %(checked) 
of all
primes~$p$.  Farmer and James \cite{FJ} have used a version of
the above result to
show that
$T_{p,k}(x)$ is
irreducible and has full Galois group
if $p<2000$ and $k\le 2000$.

In section~2 we describe the factorizations 
of $T_{p,k}^{N,\chi}(x)\bmod\ell$.  In section 3 we give
some examples.
In section 4 we  deduce some corollaries
and prove Theorem~1.
%  This is illustrated with the factorization
%of $T_{2,k}(x)\bmod 13$. 

Throughout the paper, $p$, $q$, and $\ell$ are distinct
rational primes; $N$ and $n$  are positive rational integers with
$N$, $np$, $q$, and $\ell$ pairwise coprime;  
$k$ is a positive rational integer,
and $\chi$ is a character of conductor $N$ with $\chi(-1)=(-1)^k$.  
We write $T_{n,k}^{N,\chi}(x)$ for the characteristic polynomial
of $T_n$ acting on the space of cusp forms $S_k(N,\chi)$
and put $T_{n,k}(x)=T_{n,k}^{1,1}(x)$.
When we refer to the irreducibility
or to the Galois group of a polynomial, we will mean over~$\Q$.
For prime $p$ we write $p^A||b$ to mean $p^A|b$ and
$p^{A+1}\nmid b$.

The authors thank J--P.~Serre for helpful comments concerning
the calculations in this paper.

\head
2.  Patterns in the factorization of $\TpkNchi \mod \ell$
\endhead

Fix $p$, $\ell$, $N$, and $\chi$.
In this section we describe the factorization of
$\TpkNchi \mod \ell$
as $k$ varies.  We show that those factorizations follow
a pattern.  In principle one can use this pattern to determine
$\TpkNchi \mod \ell$ for any $k$, after an initial calculation.
In the following section we give some examples
for $N=1$ and small $\ell$.

The idea here is to consider $S_k(N,\chi)\bmod\ell$ and $T_n\bmod\ell$,
and to find a basis $B_k(N,\chi)$ of $S_k(N,\chi)$
with the following nice property.

\proclaim{Lemma 1} Let $[T_n]_{k,N,\chi}$ be the matrix of $T_n$
with respect to the basis $B_k(N,\chi)$.
Then $[T_n]_{k,N,\chi}$ is block upper--triangular and
$[T_n]_{k,N,\chi}\subset [T_n]_{k+\ell-1,N,\chi}$,
where the smaller matrix is a block in the upper
left corner of the larger matrix.
In particular,
$  T_{n,k+\ell-1}^{N,\chi}(x)
\equiv g(x)\,T_{n,k}^{N,\chi}\mod \ell$
for some polynomial $g(x)$.
\endproclaim

\demo{Proof} We define the basis $B_k(N,\chi)$
as follows.
The classical congruence for the level~1 Eisenstein series 
$E_{\ell-1}\equiv 1\bmod\ell$ gives                                             
an inclusion $S_k(N,\chi)\subset S_{k+\ell-1}(N,\chi)\mod\ell$,
given by multiplication by $E_{\ell-1}$.
Choose $B_k(N,\chi)$ to respect this inclusion.
The Hecke operator $T_n$ mod~$\ell$ also respects this 
inclusion because we have  
$$
\(T_nf\)(z)={1\over n}\sum_{ad=n}\chi(a)a^{k}
\sum_{0\le b<d}f\left({az+b\over d}\right)
$$  
and $a^k\equiv a^{k+\ell-1}\mod\ell$.  This proves Lemma~1.
\enddemo

By Lemma 1, there is a sequence of polynomials 
$f_j\in (\Z/\ell\Z)[x]$, depending only on $N$, $\chi$,
$p$, $\ell$, and
($k\bmod\ell-1$), so that
$$
\TpkNchi\equiv \prod_{j=1}^J f_j(x)\bmod \ell.
\eqno(2.1)
$$
The following Proposition explains why it is easy to understand 
$\TpkNchi\mod \ell$ for any $k$.

\proclaim{Proposition 1} The sequence $f_j$ in (2.1) is periodic.
\endproclaim

Thus, the calculation of $\TpkNchi\mod\ell$ for any $k$
reduces to a finite calculation.
The periodicity of $f_j$ follows from the isomorphism
$W_{k+\ell+1}=W_k[1]$.  Here 
$W_k={\widetilde M}_{k+\ell-1}/{\widetilde M}_{k}$,
with ${\widetilde M}_{k}$ the $\Bbb F_\ell$ vector space
obtained by reducing $M_k\bmod \ell$.  
See~\cite{J}. 
We give a
different proof based on the trace formula.

\demo{Proof}  The degree of each $f_j$ is bounded by
$$
M=\max_k \(\dim S_{k+\ell-1}(N,\chi) - 
\dim S_k(N,\chi)\).
$$
This is finite because $\dim S_k(N,\chi)$ grows linearly
as a function of~$k$.

The polynomial $f_j$ has at most $M$ roots, and the coefficients
of $f_j$ are symmetric functions in those roots.
These can be expressed as polynomials in the
traces of $T_p$, $T_p^2$,..., $T_p^M$,  which in turn can be
expressed as polynomials in
$T_{p^j}$ with $j\le M$.  
Thus, the coefficients of $f_j$
are just  polynomials in the traces of 
$T_p$, $T_{p^2}$,..., $T_{p^M}\mod \ell$.

Let $\sigma_k^{N,\chi}(T_n)$ denote the trace of
$T_n$ acting on $S_k(N,\chi)$.
We need only show that 
$\sigma_k^{N,\chi}(T_n)\mod\ell\,$ is periodic
as a function of~$k$.
This follows from the Eichler--Selberg trace formula.
We state the trace formula in the following way,  
retaining only the features necessary for our proof.
For details see \cite{SV} or \cite{C}.

\proclaim{Lemma 2 (The trace formula)}  There are explicit algebraic
numbers $A$, $B_m$, and $C_d$,  depending only on
$N$, $n$, and $\chi$, 
such that
$$
\sigma_k^{N,\chi}(T_n)=
A\,n^{k/2} (k-1) \chi(\sqrt{n})
+
\sum_{|m|\le 2\sqrt{n}} B_m\,
{\eta_m^{k-1}-{\bar{\eta}}_m^{k-1}\over \eta_m-\bar{\eta}_m}
\,+
\sum_{d|n\atop 0<d\le \sqrt{n}} C_d\,d^k
,
$$
where $\displaystyle{\eta_m={m+\sqrt{m^2-4n}\over 2}}$.
Furthermore, each of
$A$, $B_m$, and $C_d$
is integral modulo the
rational prime $\ell$ for $\ell\ge 5$ and 
$\ell\nmid nN$.
\endproclaim

Each term in Lemma~2 is periodic mod~$\ell$
as a function of $k$; therefore so is
$\sigma_k^{N,\chi}(T_n)$.  
One slight complication is that
if $n$ is a square mod~$\ell$ and $\ell<4n$, 
then for some $m$
we may have $\eta_m-\bar{\eta}_m\equiv 0\mod\ell^K$ for
$K\ge 1$.  In this case choose $L$ so that
$\eta_m^{k+L}\equiv \eta_m^k\mod\ell^{K+1}$.  This proves
Proposition~1.
\enddemo

In many cases the period of $f_j$ is too large 
to actually compute on available computers.  
However, for $N=1$ and $\ell$ reasonably small, the
computation is tractable.  In the next section we
give several examples, and in the following section we
use those examples to deduce Theorem~1.

\head
3.  Some factorizations of level 1 Hecke polynomials
\endhead

We give some examples of the
factorizations of Hecke polynomials~mod~$\ell$ described in 
the previous section.  The easiest examples to calculate are for
$N=1$ and $\ell\le 7$ or $\ell=13$.  In those cases
$$
\dim S_{k+\ell-1}(1,1)\le 1+\dim S_k(1,1),
$$
so $T_{p,k}(x)\mod\ell$ factors completely
and we only need the trace of $T_p$ to determine each factor.

There is a further simplification for $\ell\le 7$ arising from the 
classification of cusp forms mod~$\ell$.
For normalized Hecke eigenforms 
$f=\sum a_nq^n$ and $g=\sum b_n q^n$, write
$f\equiv g\bmod \ell$ to mean $a_n\equiv b_n\bmod\ell$
for $(n,\ell)=1$.  For any particular $\ell$ there are
only finitely many congruence classes of Hecke eigenforms
mod~$\ell$.  These are given explicitly by Serre \cite{Ser} for 
$\ell\le 23$.

For $\ell=2$ we have $a_p\equiv 0\bmod 2$, so
$T_{p,k}(x)\equiv x^{d_k}\bmod 2$ for all $p\not=2$.
For $\ell=3$, $5$, or $7$, we have
$a_p\equiv p^m+p^n \bmod\ell$ for $p\not=\ell$, where $m$ and $n$
depend only on $\ell$ and the Hecke eigenform.
Thus, if $p\equiv q\bmod\ell$ then
$a_p\equiv a_q\bmod\ell$,
so $T_{p,k}(x)\equiv T_{q,k}(x)\bmod\ell$ if $p\equiv q\bmod\ell$.
For $\ell\le 7$ this reduces the determination of
$T_{p,k}(x)\bmod\ell$ for all $p$ and $k$ to the calculation of a few
cases.  The results of those calculations are presented in the
following Theorem.

\proclaim{Theorem 2a}  For prime $\ell\le 7$ and $p\not =\ell$,
the Hecke polynomial
$T_{p,k}(x)$ factors as
$$
T_{p,k}(x)\equiv \prod_{j=1}^{d_k} (x-a_j) \bmod \ell
$$
where $\{a_j\}$ is a periodic sequence depending only on 
$p\bmod\ell$ and $k\bmod(\ell-1)$.
We have
$$
T_{p,k}(x)\equiv \cases (x-2)^{d_k}\bmod 3, &\hbox{if } p\equiv 1\bmod 3\\
x^{d_k}\bmod 3& \hbox{if }p\equiv 2\bmod 3\\
\endcases
$$
For $\ell=5$ or $7$ the results are summarized in the tables below.
The rows are labeled by the smallest prime in each congruence
class mod~$\ell$, the columns are labeled by the congruence
class of $k\bmod (\ell-1)$, and the table entry gives one
period of the sequence $\{a_j\}$.
The tables are sufficient to determine $T_{p,k}(x)\bmod\ell$ for
all $p$ because if $\ell\le 7$ and $p\equiv q\bmod\ell$ then
$T_{p,k}(x)\equiv T_{q,k}(x)\bmod\ell$.
\endproclaim
\hbox{\rm %\smallertype
\vbox{\tabskip=0pt \offinterlineskip
\halign to \tablewidthA{
\newstrut
#& #\tabskip=0em plus 1em&
  \hfil#& \vrule#& #\hfil& \newstrutA#&
  #\hfil&\hfil#&
#\tabskip=0pt\cr     %\tablerule
\makeheadingfortableA
&&$p=11$&&(2)&&(2)& \cr
&&2&&(1,4)&&(2,3)& \cr
&&3&&(2,3)&&(1,4)& \cr
&&19&&(0)&&(0)&\cr
&\newstrutB\cr
&\newstrutB\cr
&\newstrutB\cr
}
}
\hskip .2in
\vbox{\tabskip=0pt \offinterlineskip
\halign to \tablewidthB{
\newstrut
#& #\tabskip=1em plus 2em&
  \hfil#& \vrule#& #\hfil& \vrule#&
  #\hfil& \newstrutA#&
  #\hfil&\hfil#&
#\tabskip=0pt\cr     %\tablerule
\makeheadingfortableB
&&$p=29$&&(2)&&(2)&&(2) \cr
&&2&&(4,5)&&(1,3)&&(6,2) \cr
&&3&&(0,1,0,6)&&(0,3,0,4)&&(5,0,2,0)\cr
&&11&&(1,3)&&(4,5)&&(6,2)\cr
&&5&&(0,3,0,4)&&(0,1,0,6)&&(2,0,5,0)\cr
&&13&&(0)&&(0)&&(0)\cr
&\newstrutB\cr
}}
}

%To illustrate the factorizations in a different way, we have
%$$\eqalignno{
%T_{2,k}(x)&\equiv (x-1)^a(x-4)^b \bmod 5 &\hbox{if } k\equiv0\bmod 4 \cr
%\cr
%T_{3,k}(x)&\equiv x^m(x-3)^a(x-4)^b \bmod 7 &\hbox{if } k\equiv2\bmod 6,\cr
%}$$

For $\ell>7$ there are no simple relationships between
$T_{p,k}$ and $T_{q,k}\mod\ell$, 
so  each $T_{p,k}\mod\ell$  must be calculated separately.
We give an example mod~13.

\proclaim{Theorem 2b}  We have a factorization                                  
$T_{2,k}(x)\equiv \prod_{j=1}^{d_k}(x-a_j)\bmod 13$, where                      
$\{a_j\}$ is a sequence of period~14 depending only                             
on $k\bmod 12$.   The first 14 terms of each sequence are given                 
in the following table, where each row corresponds to                           
a congruence class of $k\bmod 12$.  \hfill\break
%\endproclaim   
\hbox{\rm
\hskip .5in
\vbox{\tabskip=0pt \offinterlineskip   
\vskip .2in
\halign to \tablewidthC{                                                        
\newstrut                                                                       
#& #\tabskip=1em plus 2em&                                                      
  \hfil#& \newstrutA#& #\hfil& #&                                               
%  #\hfil& #&                                                                   
%  #\hfil&\hfil#&                                                               
#
\tabskip=0pt\cr     %\tablerule                                                 
\makeheadingfortableC                                                           
&&$k\equiv0\bmod12$&&                                                           
(2, 12, 9, 4, 1, 11, 5, 11, 1, 4, 9, 12, 2, 8)\cr                               
&&2\phantom{\hbox{ mod 12}}&&                                                   
(4, 11, 5, 8, 2, 9, 10, 9, 2, 8, 5, 11, 4, 3)\cr                                
&&4\phantom{\hbox{ mod 12}}&&                                                   
(8, 6, 8, 9, 10, 3, 4, 5, 7, 5, 4, 3, 10, 9)\cr                                 
&&6\phantom{\hbox{ mod 12}}&&                                                   
(5, 3, 12, 3, 5, 7, 6, 8, 10, 1, 10, 8, 6, 7)\cr                                
&&8\phantom{\hbox{ mod 12}}&&                                                   
(1, 10, 6, 11, 6, 10, 1, 12, 3, 7, 2, 7, 3, 12)\cr                              
&&10\phantom{\hbox{ mod 2}}&&                                                   
(11, 2, 7, 12, 9, 12, 7, 2, 11, 6, 1, 4, 1, 6)                                  
\cr 
&\newstrutB\cr                                                                  
}}   
}
\endproclaim

The calculations for Theorem~2 were done in Mathematica.
The method was to find a basis for $S_k(1)$
using $\Delta$, $E_4$, and  $E_6$,
and then explicitly compute the action of the Hecke operator.
Representative cases were checked using the trace formula.

The period of $a_j$ given in Theorem~2 is shorter than might 
have been predicted from Lemma~2.  
If $p$ is not a square mod~$\ell$
then $(m^2-4p,\ell)=1$ for all $m$, so Fermat's little theorem 
gives
$\eta_m^k\equiv\eta_m^{k+\ell^2-1}\mod\ell$,
which implies
$\sigma_k(T_p)=\sigma_{k+\ell^2-1}(T_p)\bmod \ell$.
This implies that
the $a_j\bmod\ell$ has period at most~$(\ell^2-1)/12$.
This is the actual period given for those cases in Theorem~2a
and Theorem~2b.
If $p$ is a square
mod~$\ell$ then the period of $\sigma_k(T_p)$ is much larger,
leading to a larger upper--bound for the period of $a_j$.
(The largest case we needed is $p=11$, $m=3$, $\ell=7$, 
where we have $\eta_3^{295}\equiv \eta_3\bmod 7^2$.)
However, as can be seen in the examples, the  period
of $a_j$ is actually smaller.  
J--P.~Serre has pointed out to us that the elementary
result $W_{k+p+1}=W_k[1]$ gives the indicated periodicity,
and this can also be used to explain the various patterns
which appear in the table in Theorem~2b.

\head
4.  Proof of Theorem 1
\endhead

In this section we deduce some consequences of the factorizations
given in Theorem~2 and we use those factorizations to prove Theorem~1.
%We use those factorizations to deduce the following corollary.

\proclaim{Corollary}  Suppose for some $n$ that
$T_{n,k}(x)$ is irreducible over~$\Q$.  Then $T_{p,k}(x)$
is irreducible if either of the following holds:
\itemitem{i)} $\dim S_k(1)$ is odd and
$p\not\equiv \pm 1\bmod 5$
or
$p\not\equiv \pm 1\bmod 7$, or
\itemitem{ii)} $\dim S_k(1) \equiv 2\bmod 4$ and 
$p\equiv 3$ or $5\bmod 7$.
\endproclaim

Note that the Corollary applies to $1/2$ of all pairs $(k,p)$.
Both the Corollary and Theorem~1 follow from the factorizations 
in Theorem~2 and the following
proposition.

\proclaim{Proposition 2}  Suppose $T_{n,k}(x)$ is irreducible
for some $n$.  Then for each $m$
we have $T_{m,k}(x)=f(x)^r$ with $f(x)$ irreducible
and $r\in\N$.
Suppose  
$T_{n,k}(x)$ is irreducible and has full 
Galois group for some~$n$.
Then for each $m$
either $T_{m,k}(x)$ is irreducible and has full Galois group,
or $T_{m,k}(x)=(x-a)^{d_k}$ for some $a\in \Z$.
\endproclaim

\demo
{Proof of the Corollary and Theorem~1}
If $\,T_{p,k}(x)=f(x)^r$ then each root of $T_{p,k}(x)\bmod\ell$
has multiplicity divisible by $r$.  Using $\ell=5$ or~7
and the factorizations in Theorem~2a  gives $r=1$ for all cases
covered by the Corollary.

If $T_{p,k}(x)=(x-a)^r$ then $T_{p,k}(x)$ has only one root
mod~$\ell$.  By Theorem~2a this 
does not hold in the cases covered by Theorem~1.
\enddemo

\demo
{Sketch of Proof of Proposition~2}
The idea is to consider how the Galois group of the field of Fourier
coefficients of the cusp forms acts on the Hecke basis.
 
Let $B_k$ be the set of normalized Hecke eigenforms in $S_k(1)$,
let $\K/\Q$ be the field generated by the set of Fourier
coefficients of the $f\in B_k$, and let $G=Gal(\K/\Q)$.
If $f(z)=\sum a_nq^n$ with $\{a_n\}\subset \K$ and $\sigma\in G$
then set $\sigma f(z)=\sum \sigma a_nq^n$. 

Since $S_k(1)$
has a rational basis, $G$ acts on $S_k(1)$.
The group $G$ also acts on the set $B_k$ because $G$ commutes
with the Hecke operators.   This induces an action on the
roots of $T_{n,k}(x)$ because those roots are the $n$th Fourier
of the Hecke eigenforms.  That action is the same as the natural
action given by the inclusion of those roots in $\K$.
Thus, if $G$ acts transitively on $B_k$ then $G$ acts 
transitively on the roots of $T_{m,k}(x)$, which implies the
first conclusion in the Proposition.
If $G$ acts as the full symmetric group on $B_k$, then $G$ acts
as the full symmetric group on the roots of $T_{m,k}(x)$,
which implies the second conclusion in the Proposition.

To finish the proof, note that
the Galois group of $T_{n,k}(x)$ acts on $B_k$ by extending
the action on the $n$th Fourier coefficients.  Thus, if
$T_{n,k}(x)$ is irreducible then $\K$ is the splitting
field of $T_{n,k}(x)$ and $G$ is its Galois group. 
This proves Proposition~2.
\enddemo

It would be interesting to prove a version of Maeda's
conjecture of the form ``if $T_{n,k}(x)$ is irreducible for some $n$,
then it is irreducible for all $n$,'' 
or even ``if $T_{n,k}(x)$ is irreducible for some $n$,
then $T_{2,k}(x)$ is irreducible.''
The methods
used here are only able to establish such a result for a large
proportion of weights~$k$.
By combining Theorem~2b with the Corollary above, we have that
if $T_{n,k}(x)$ is irreducible for some $n$ and
$\dim S_k(1)$ is not a multiple of 14, then
$T_{2,k}(x)$ is irreducible, and if 
$\dim S_k(1)$ is not a multiple of 28, then
$T_{2,k}(x)$ or $T_{3,k}(x)$ is irreducible.
Establishing factorizations of $T_{2,k}(x) \mod\ell$ for
larger $\ell$ would increase the proportion of $k$ for
which such a result is known.  This is 
reminiscent of Lehmer's conjecture that the
Ramanujan $\tau$--function never vanishes, where the results
in that direction come from congruence properties
of $\tau(n)$.

\Refs

\item{[B]} {\sl K.~Buzzard},  On the eigenvalues of the Hecke operator
   $T\sb 2$, J.~Number Theory  {\bf 57} (1996), no. 1, 130--132.

\item{[C]} {\sl H.~Cohen},  Trace des op\'erateurs de Hecke
sur $\Gamma_0(N)$, in ``S\'eminaire de Th\'eorie des
nombres,'' exp. no.~4, Bordeaux, 1976-1977.

\item{[CF]} {\sl J.B.~Conrey} and {\sl D.W.~Farmer},  Hecke
operators and the nonvanishing of $L$--functions, to appear in
the proceedings of the Penn State conference, Summer 1997.

\item{[FJ]} {\sl D.W.~Farmer} and {\sl K.~James}, The irreducibility
of some level~1 Hecke polynomials, preprint.

\item{[HM]} {\sl H.~Hida} and {\sl Y.~Maeda}, Non--abelian base
change for totally real fields, Pac.~J.~Math {\bf 181} no.~3
(1997) 189--218.

\item{[JO]} {\sl K.~James} and {\sl K.~Ono}, 
On the factorization of Hecke polynomials,
preprint.

\item{[J]} {\sl N.~Jochnowitz}, Congruences between systems of
eigenvalues of modular forms, Trans. Amer. Math. Soc. {\bf 270}
(1982), 269--285.

\item{[KZ]} {\sl W.~Kohnen} and {\sl D.~Zagier},
Values of $L$-series of modular forms at the center of the critical strip,
Invent. Math. {\bf 64} (1981), no. 2, 175--198.

\item{[Ser]} {\sl J--P Serre}, Valeurs propres des op\'erateurs
de Hecke modulo~$\ell$.  Ast\'erisque {\bf 24}--{\bf 25} (1975), 109--117.

\item{[SV]} {\sl R.~Schoof} and {\sl M.~van der Vlugt},
Hecke operators and the weight distribution of certain codes,
Journal of Combinatorial Theory A, {\bf 75} (1991), 163--186.

\endRefs

\bye
\enddocument